\documentstyle[aps,preprint,overcite,graphicx,epsf]{revtex}

\tightenlines

\def\affiliation{\address}
\def\wideabs#1{#1}
\newcommand{\bfi}{\bfseries\itshape}
\def\bfnabla{\mbox{\boldmath $\nabla$}}
\def\bfeta{\mbox{\boldmath $\eta$}}
\def\be{\begin{equation}}
\def\eeq{\end{equation}}
\def\bea{\begin{eqnarray}}
\def\eea{\end{eqnarray}}

\def\k1{k_1}

\def\a{\alpha}
\def\fa{\left(1-\a^2\triangle\right)}

\def\v#1{\text{\bfi #1}}

\begin{document}

\wideabs{
\title{Enhancement of the inverse-cascade of energy in the two-dimensional averaged Euler equations}
\author{B.T. Nadiga}
\affiliation{Earth and Environmental Sciences, MS-B296 \\
	Los Alamos National Laboratory, Los Alamos, NM 87545}
\author{S. Shkoller}
\affiliation{Department of Mathematics \\
University of California, Davis CA 95616}
\maketitle
\begin{abstract}
\hspace{-.6em} 
For a particular choice of the smoothing kernel, it is shown that
the system of partial differential equations governing the vortex-blob method
corresponds to the averaged Euler equations. These latter equations
have recently been derived by averaging the Euler equations
over Lagrangian fluctuations of length scale $\a$, and the same system is
also encountered in the description of inviscid and incompressible flow 
of second-grade polymeric (non-Newtonian) fluids.  While previous
studies of this system have noted the suppression of nonlinear interaction
between modes smaller than $\a$, we show that
the modification of the nonlinear advection term also 
acts to enhance the inverse-cascade of energy in two-dimensional turbulence
and thereby affects scales of motion larger than $\a$ as well. This latter
effect is reminiscent of the drag-reduction that occurs in a turbulent
flow when a dilute polymer is added.
\end{abstract}
}

\newpage
The two-dimensional incompressible, Euler equations are
\be\label{01}
\partial_t \omega + \bfnabla\cdot\left(\v{u}\,\omega\right) =0, \quad
\bfnabla\cdot\v{u}=0, \quad\omega(t=0)=\omega_0,
\eeq
where $\omega=\bfnabla\times\v{u}$ 
is the vorticity, $\v{u}$ is the spatial velocity
vector field, $t$ denotes time, and all the dependent variables depend on $t$
and $\v{x}=(x_1,x_2)$, the
Cartesian coordinates in the plane.  An inversion of the 
vorticity-velocity relation yields
$\v{u} = \int \v{K}(\v{x},\v{y}) \omega(\v{y})d\v{y}$,
where $\v{K}= \bfnabla^\perp G$, $G$ is the solution of $-\triangle G= \delta$,
and $\bfnabla^\perp= (-\partial_{x_1}, \partial_{x_2})$.For 
fluid motion over the
entire plane, $\v{K}(\v{x},\v{y})=(2\pi)^{-1}\bfnabla^\perp\log|\v{x}-\v{y}|$.
Let $\bfeta_t$ denote the flow of $\v{u}_t=\v{u}(t,\cdot)$, so that
$d\bfeta_t/d t = \v{u}_t(\bfeta(t))$.
Because $\v{u}_t$ is divergence-free, the flow map $\bfeta_t$ is an 
area-preserving transformation for each $t$.  It follows that
\begin{eqnarray}\label{04}
\frac{d\bfeta_t}{d t} & = \int \v{K}(\bfeta_t(\v{x}),\bfeta_t(\v{y})) 
\omega(\bfeta_t(\v{y}))d\v{y}
                          \nonumber\\
                    & = \int \v{K}(\bfeta_t(\v{x}),\bfeta_t(\v{y})) 
             \omega_0(\v{y})d\v{y},
\end{eqnarray}
where the last equality is a consequence of the 
pointwise conservation of vorticity
along Lagrangian trajectories, $\omega(\bfeta_t(\v{x}))= \omega_0(\v{x})$.
Thus, the initial vorticity field completely determines the
fluid motion.  Choosing the initial vorticity to be a sum of $N$ point vortices
$\delta_i$ positioned at the points $\v{x}_i$ in the plane with circulations
$\Gamma_i$, $\omega_0=\sum_{i=1}^N \Gamma_i \delta_i$, equation (\ref{04})
produces the classical point-vortex approximation to (\ref{01}).  This
approximation is known to be highly unstable, as finite-time collapse
of vortex centers may occur~\cite{MP}.

Chorin's vortex blob method~\cite{C} alleviates the instability of the
point-vortex scheme by smoothing each delta function $\delta_i$ with a vortex
blob $\chi$, a function that decays at infinity, and whose mass is
mostly supported in a disc of diameter $\alpha$.  Thus, instead of using the
integral kernel $\v{K}(\v{x},\v{y})$, one uses the smoother kernel $\v{K}^\alpha
= \bfnabla^\perp G^\alpha$ where $G^\alpha$ is the solution of
$ -\triangle G^\alpha = \chi.$
The vortex-blob method then evolves the point-vortex initial data, which we
shall now call $q_0$, by the ordinary differential equation
\begin{equation}\label{05}
\frac{d\bfeta_t^\alpha}{d t} = \int \v{K}^\alpha(\bfeta_t^\alpha(\v{x}),
\bfeta_t^\alpha(\v{y}))
q_0(\v{y})d\v{y}.
\end{equation}
Henceforth, to keep the notation concise, we will drop the superscript 
$\a$ when there is no ambiguity.

When the vortex-blob $G^\alpha$ is the modified Bessel function of the
second kind $K_0$, it is the fundamental solution of the operator 
$(1-\alpha^2 \triangle)$ in the plane, 
and the vorticity $q$ is related to the
smoothed velocity vector field $\v{u}$ by $q=(1-\alpha^2 \triangle)
\bfnabla\times\v{u}$.  
Thus, Chorin's vortex method for this choice of
smoothing is given by the partial differential equation
\be\label{06}
\partial_t q + \bfnabla\cdot\left(\v{u}\, q\right)=0,
\quad\bfnabla\cdot\v{u}=0, \quad q(t=0)=q_0,
\eeq

The system of equations (\ref{06}) is also known as the two-dimensional
isotropic
averaged Euler equations, and are derived by averaging over Lagrangian
fluctuations of order $\alpha$ about the macroscopic flow 
field~\cite{HMR,Shk1,MS}.
When the constant $\alpha>0$ is interpreted as a material
parameter which measures the elastic response of the fluid due to 
polymerization instead of as a spatial length scale, then (\ref{06})
are also exactly the equations that govern the inviscid flow
of a second-grade non-Newtonian fluid~\cite{NT}. 
According to Noll's theory of simple materials, 
(\ref{06}) are obtained from the
unique constitutive law that satisfies material frame-indifference and 
observer objectivity.  Consequently, the vortex method with the 
Bessel function $K_0$ smoothing naturally inherits these  
characteristics~\cite{OS}.

We find the connections between averaging Euler equations over Lagrangian
fluctuations, a constitutive theory for polymeric fluids,
and a classical numerical algorithm to be quite
intriguing and suspect that these equations will be important 
from a modeling standpoint.
However, most previous studies of the averaged Euler equations have been of
a mathematical nature, and we are aware of only a few cases where 
this system has been used as a (dynamic) modeling tool: Chen et al.~\cite{Chen}
used a viscous version of the three-dimensional averaged Euler equations
to simulate isotropic turbulence and found that they could reproduce
large-scale features without fully resolving the flow.
Nadiga~\cite{N} considered the inviscid two-dimensional form of the
averaged Euler equations and demonstrated that for suitably chosen values
of $\a$, the large-scale {\em spectral-scalings} of the Euler equations
could be preserved while achieving a faster spectral decay at the 
smaller scales. Finally, Nadiga \& Margolin~\cite{NM} used an extension of
the two-dimensional averaged Euler equations in a geophysical
context to model the effects of mesoscale eddies on mean flow.

Before we go on to consider numerical simulation of the averaged Euler system,
we wish to point out that there is also a beautiful geometric
structure to (\ref{06}) which follows the framework developed
by Arnold~\cite{A} and Ebin and Marsden~\cite{EM}.  
While the details of this particular issue are
far outside the scope of this article, it is, nevertheless, worthwhile to state
the result.  Arnold showed that the appropriate configuration space for a
perfect incompressible fluid is the group of all area preserving diffeomorphisms
of the fluid container, and that solutions of the Euler equations are geodesics
on this group with respect to a certain kinetic energy metric, characterized
by the inner-product $\int \left(\v{u}\cdot\v{v}\right)\, d\v{x}$ 
for two divergence-free vector fields $\v{u}$ and $\v{v}$.  
The system (\ref{06}) also has this geometric property, but
now the metric is instead characterized by 
$\int \left(\v{u}\cdot\v{v} + 2 \alpha^2 \text{Def}(\v{u}) \cdot 
\text{Def}(\v{v})\right) d\v{x},$
where $\text{Def}(\v{u})$ is the rate of deformation tensor 
$(\bfnabla\v{u} + (\bfnabla\v{u})^T)/2$~\cite{Shk2}.
Equations (\ref{06}) thus preserve the Hamiltonian structure
of the Euler equations.  In particular, vorticity remains pointwise conserved
by the smooth Lagrangian flow $\bfeta_t^\alpha$ so that
$ q(\bfeta_t^\alpha(\v{x}))= q_0(\v{x}),$
the vorticity momenta
$\int q^p d\v{x}$
are conserved, and so the Kelvin circulation theorem remains 
intact~\cite{HMR} as well. 

Since, in each of the three different scenarios---averaged Euler,
vortex-blob method, and inviscid second-grade fluid---the
essential modification of the original equations
is a change of the advective nonlinearity
of fluid dynamics, we will now consider 
forced-dissipative simulations of the system (\ref{06}) to demonstrate
the effect of such an inviscid modification.
If we use a vorticity-stream function formulation, 
the evolution of both the Euler and averaged Euler sytems
can be represented by
\be\label{general}
\frac{\partial\omega}{\partial t} + (1-\alpha^2 \triangle)^{-1} J[\psi, 
(1-\alpha^2 \triangle) \omega] = F + D;
\quad 
\end{equation}
where $\omega=\triangle\psi$, $J$ is the Jacobian, $F$ the forcing, and
$D$ the dissipation.
(The Euler system corresponds to $\a=0$.)
Our numerical scheme consists of 
a fully dealiased pseudospectral spatial discretization and
a (nominally) fifth-order, adaptive timestep, embedded
Runge-Kutta Cash-Karp temporal discretization of (\ref{general})
(see Nadiga~\cite{N} for details).
With such a scheme, among the infinity of inviscid ($F=D=0$) 
conserved quantities for (\ref{general}), 
the only two conservation properties that survive are those for
the kinetic energy $E_{H^1}$ and enstrophy  $Z_{H^2}$ given respectively by
\bea\label{elrke}
E_{H^1}&=&\frac{1}{2}\int\left(|\v u|^2 + \alpha^2 |\bfnabla \v u|^2\right)
d\v{x}
\;\left(=\|\v u\|_{H^1}^2\right),\nonumber\\
Z_{H^2}&=& \frac{1}{2} \int \left[\fa\omega\right]^2 d\v{x}
\;\left(=\|\omega\|^2_{H^2}\right).
\eea

In the forced-dissipative runs to be considered,
the forcing $F$ is achieved by keeping the amplitudes of modes with wavenumbers
in the small wavenumber band
$10\leq k < 10.001$ constant in time. The dissipation, $D$, is a combination
of a fourth order hyperviscous operator and a large-scale friction term:
$D = \delta\psi-\left(-\nu\triangle\right)^4\omega$,
as has been used in numerous previous studies of two-dimensional turbulence.
The form and value of the forcing and dissipation are held exactly the same
for all the runs to be presented, irrespective of the resolution and the 
value of $\a$.

\begin{figure}
\epsfxsize=\columnwidth
\centerline{\epsfbox{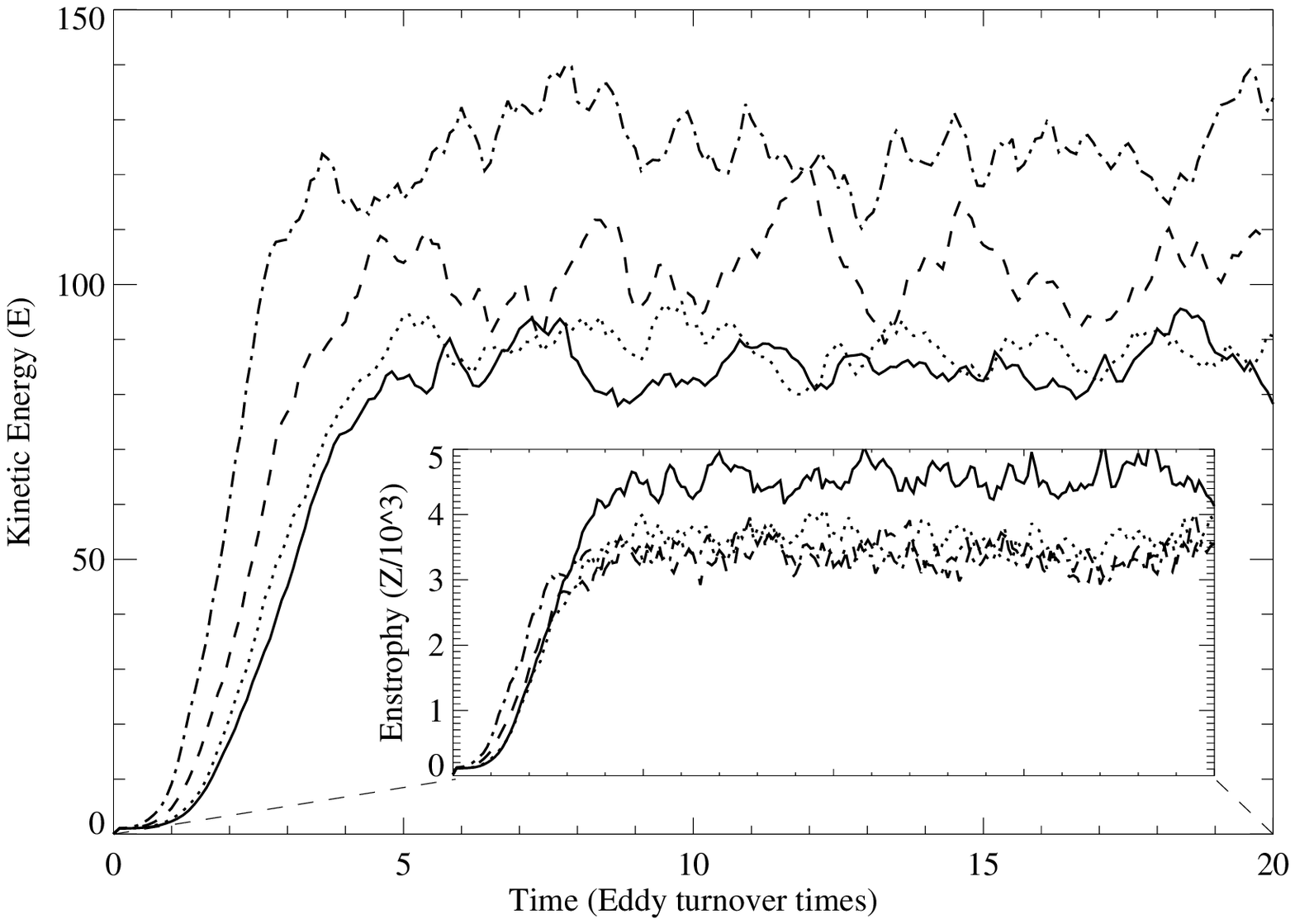}}
\caption{The evolution of kinetic energy E with time for
$k_\a=\infty$ (solid line), $k_\a=42$ (dotted line), $k_\a=21$ (dashed line), 
and $k_\a=14$ (dot-dashed line). An increase in $\a$, for identical
forcing and dissipation, results in an overall reduced viscous behavior.
In the inset is shown the evolution of enstrophy Z for the same time interval
and for the same four values of $\a$ and with the same line types as for
kinetic energy.
While there is a  significant difference
between zero and non-zero $\a$ cases, the dependence on the actual
value of $\a$ itself is rather weak.}
\label{ul2}
\end{figure}

On the one hand, it could be argued that 
since the energy and enstrophy 
that are conserved (in an unforced-inviscid setting) 
are $E_{H^1}$ and $Z_{H^2}$ respectively, it is their dynamics 
which is of primary importance.
On the other, it could be argued that in the context of (\ref{general}), 
the interest in small scales is only in so much as it affects the 
larger scales and to that extent $\a$ has no primary significance and
that it is really the large scale components of energy and
enstrophy in (\ref{elrke}) that are of primary interest.
While both these points of view are reasonable, in this 
short article, we proceed with the latter
and concern ourselves
with the dynamics of the usual kinetic energy and usual enstrophy  
as given by
\be\label{ke}
E=\frac{1}{2}\int|\v{u}|^2 d\v{x},\quad
Z = \frac{1}{2} \int \omega^2 d\v{x}.
\end{equation}
Fig.~\ref{ul2} shows the evolution of the 
kinetic energy $E$ with time
for four different values of $k_\a$.
For these computations, 512 physical grid points were used in 
each direction, resulting in, after accounting for dealiasing,
a maximum, circularly-symmetric wavenumber, $k_{max}$, of 170.
The four runs correspond to $k_\a$ of $\infty$ (dissipative Euler)
and 42, 21, and 14.
This figure shows 
that for identical forcing and dissipation, the tendency with increasing
$\a$ (equivalently decreasing $k_\a$) is to achieve an overall balance 
which makes the flow less viscous. 

While the kinetic energy of the runs with different $\a$ shows 
a definite trend (increasing with increasing $\a$), such is not the
case with the enstrophy shown for the same four cases in the inset 
of Fig.~\ref{ul2}.
Here interestingly, all the runs with nonzero $\a$ seem to display 
approximately the same level of enstrophy which is lower than for
$k_\a=0$.  This indicates: 
\begin{itemize}
\item That the small-scale
behavior is quite different when $\a=0$ and when $\a$ is nonzero
(as noted in  previous studies~\cite{HMR,Chen,N}), but 
that this difference is not sensitively dependent on the value of $\a$ for the
interesting range of values of $\a$. 
\item That the more significant change 
with $\a$ is the behavior of the large scales.
\end{itemize}

\begin{figure}
\epsfxsize=\columnwidth
\centerline{\epsfbox{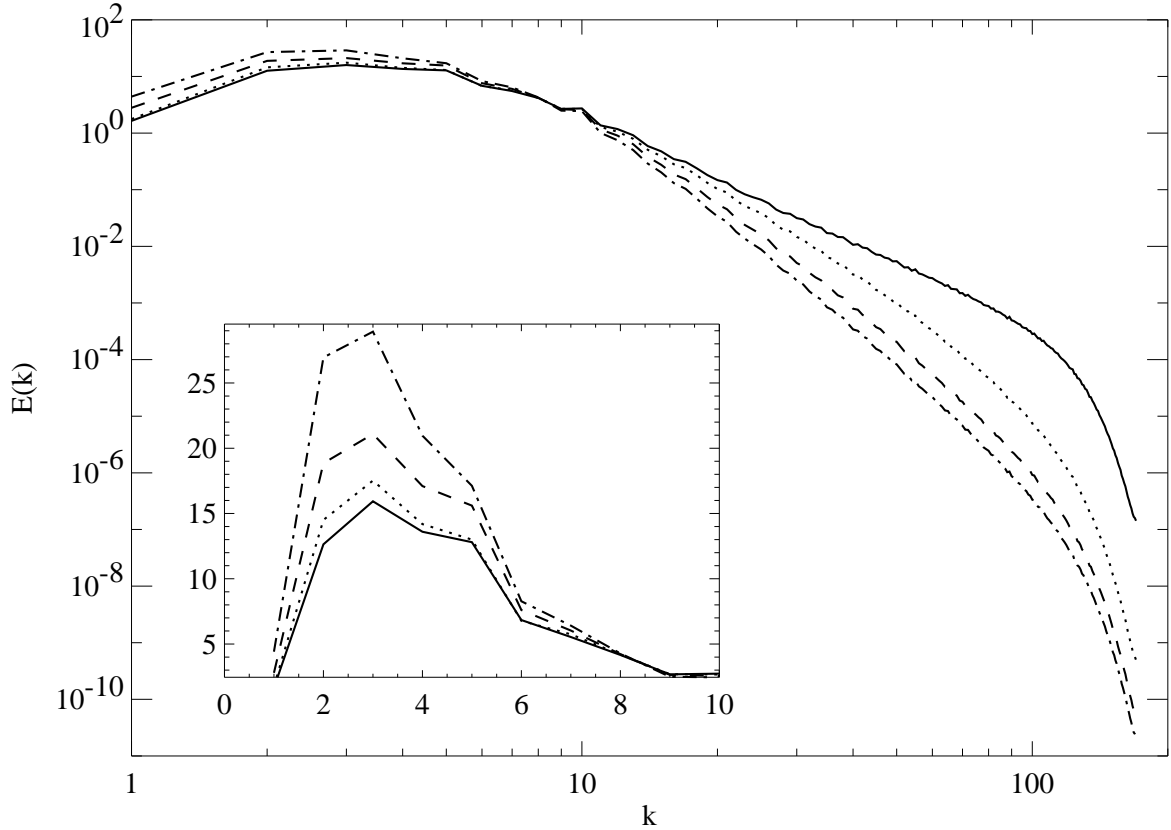}}
\caption{Stationary wavenumber-energy spectra (log-log scale) for the 
forced-dissipative simulations of the averaged Euler
equations with zero and nonzero $\a$. 
$k_\a=\infty$ (solid line), $k_\a=42$ (dotted line), $k_\a=21$ (dashed line), 
and $k_\a=14$ (dot-dashed line). The insert shows the same plot with a
linear-linear scale for the first ten wavenumbers.
The enhanced inverse cascade of energy 
and the suppressed energy level at smaller scales with increasing $\a$ is
evident.}
\label{spec1}
\end{figure}

Therefore, to further examine the nature of this
(reduced-viscous) behavior of the large scales, 
we examine the energy-wavenumber spectra 
in Fig.~\ref{spec1}. Here, the average of the 
one dimensional energy spectrum $E(k)$ between times 5 and 20 is
plotted against the scalar wavenumber $k$. Figure~\ref{spec1} shows that 
the reduced-viscous behavior for increasing $\a$ is achieved by 
systematically increasing the energy in modes larger in scale than
the forcing scale and decreasing the energy in modes smaller in scale 
compared to the forcing scale.  

The larger energy content in the
larger scales implies an enhancement of the
inverse cascade of energy of two-dimensional turbulence by the 
nonlinear-dispersive modification of the advective
nonlinearity when $\a>0$ in (\ref{general}).  So also, the decreased
energy content in the smaller scales is attributable
to the same nonlinear-dispersive modification.
In the following, we give a simple dimensional argument to explain the
observed behavior.
For this, consider the governing equations in the form (\ref{06}).
In close analogy with the classical picture for the inertial ranges of
two-dimensional dissipative Euler equations~\cite{K}, a 
Kolmogorov-like cascade picture for (\ref{06}) shows that the inertial
range consists of two subranges, the enstrophy cascade subrange where
there is a down-scale cascade of the $Z_{H^2}$ enstrophy defined 
in (\ref{elrke}), and the energy cascade subrange where there is an
up-scale cascade of the $E_{H^1}$ energy defined in (\ref{elrke}).
$E_{H^1}$ and $Z_{H^2}$ are the relevant energy and enstrophies
since these are the ones which are conserved in an inviscid and 
unforced case.

If we assume that the wavenumber 
$k_\a$ only appears in the Helmholtz operator, as it does in the
governing equations, then 
\begin{itemize}
\item In the enstrophy cascade subrange,
\be\label{k1}
E(k)\sim\beta_{H^2}^a k^b,
\end{equation}
where $\beta_{H^2}$ is the rate of dissipation of $Z_{H^2}$ enstrophy,
and $a$ and $b$ are exponents to be determined by dimensional analysis.
If $L$ and $T$ are characteristic length and time scales in the 
enstrophy cascade subrange, (\ref{k1}) implies
$$ L^3 T^{-2} = T^{-3a} \left(1+\a^2 L^{-2}\right)^{2a} L^{-b} , $$
from which $a=2/3$. However, even in the enstrophy cascade subrange,
the value of $b$ depends on the the ratio $\a/L$. For $\a\ll L$, 
of course, $b=-3$, and the classical~\cite{K} $E(k)\sim k^{-3}$
is recovered; when $\a\gg L$, $E(k)\sim k^{-\frac{17}{3}}$.
Finally, when $\a$ is comparable to $L$, it is easy to see that
$E(k)$ decays faster than for Euler,
but slower than $k^{-\frac{17}{3}}$
(as may be seen in Fig.~\ref{spec1}).
\item In the energy cascade subrange,
\be\label{k2}
E(k)\sim\epsilon_{H^1}^a k^b,
\end{equation}
where $\epsilon_{H^1}$ is the rate of dissipation of $E_{H^1}$ energy,
and $a$ and $b$ are exponents to be determined by dimensional analysis.
If $L$ and $T$ are characteristic length and time scales, now, in the 
energy cascade subrange, (\ref{k2}) implies
$$ L^3 T^{-2} = T^{-3a} \left(1+\a^2 L^{-2}\right)^a L^{2a-b} , $$
from which $a=2/3$. Again, even in the energy cascade subrange,
the value of $b$ depends on the the ratio $\a/L$. For $\a\ll L$, 
of course, $b=-\frac{5}{3}$, and the classical~\cite{K} 
$E(k)\sim k^{-\frac{5}{3}}$
is recovered; when $\a\gg L$, $E(k)\sim k^{-3}$.
When $\a$ is comparable to $L$, it is easy to see that the inverse cascade
of energy is enhanced, that is
$E(k)$ increases with decreasing $k$ faster than $k^{-\frac{5}{3}}$ (Euler) 
but slower than $k^{-3}$.
\end{itemize}

While the effects of the enhancement of the inverse cascade of energy is
clear in Fig.~\ref{spec1},
we defer the verification of the asymptotic values of the exponent $b$ to
later studies when we can afford much larger simulations with a
good dynamic range in each of the inertial subranges.

\begin{figure}
\epsfxsize=\columnwidth
\centerline{\epsfbox{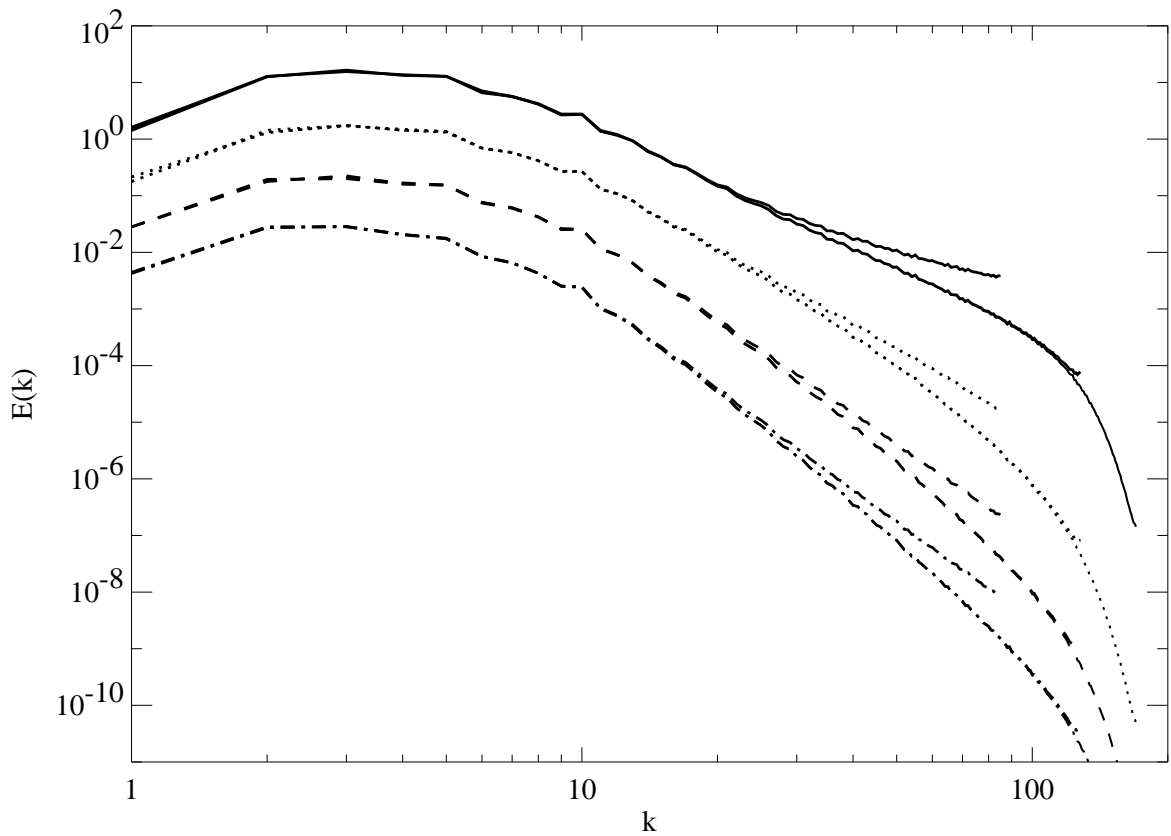}}
\caption{The spectra for the four cases in Fig.~\ref{spec1} are replotted
along with the spectra for the same four cases with resolution reduced by
25\% and 50\% in each direction. 
The sets of spectra for each $\a$ are offset by a decade each to 
improve clarity. The degree of non-resolution of the flow with the
reduced resolution is indicated by the difference between that spectrum and
the spectrum for the fully resolved case. With a 25\% reduction in resolution,
the flows are almost resolved, while with a 50\% reduction, the flows are
not fully resolved anymore.
The degree of non-resolution is independent of $\a$ to the lowest order.}
\label{spec2}
\end{figure}

The steeper fall-off of the energy spectrum with $k$ in the
enstrophy cascade range of wavenumbers when $\a>0$, compared to
Euler may, at first, suggest that a coarser resolution
may be sufficient to resolve the flow when $\a>0$ (for the
same forcing and dissipation).
However, this is not the case, as should be clear from Fig.~\ref{spec2}.
In this figure, the spectra of the cases previously discussed is
replotted together with the corresponding spectra when the
resolution is reduced by a quarter ($k_{max}=128$) and a half ($k_{max}=85$).
(The spectra for the different values of $\a$ 
are offset to improve clarity.) The degree of non-resolution
of the flow due to the reduced resolution is indicated by the 
deviation of that spectrum from that for the fully resolved case.
With a 25\% reduction in resolution,
the flows are almost resolved for all values of $\a$, 
while with a 50\% reduction, the flows are
not fully resolved anymore.
Importantly, the degree of non-resolution is independent 
of $\a$ to the lowest order.


Besides their use in describing mean motion, the averaged Euler
equations have arisen independently in at least two other
contexts---second grade polymeric fluids and vortex blob methods.
In this note, we make two observations that are likely to be of fundamental 
importance in understanding the relevance of these models in 
describing more realistic flows:
While it has been previously noted that with these equations,
nonlinear interactions at scales small compared to $\a$ are
suppressed, we have shown here that the modification of the nonlinear
advection term in these equations also leads to an enhancement of the
inverse cascade of energy in two dimentions---a characteristic feature
of two-dimensional turbulence.  This in turn implies (1) an overall
reduced-viscous behavior and (2) a significant modification of the
dynamics of scales larger than $\a$, both reminiscent of the
phenomenon of drag reduction in a turbulent flow when a dilute polymer
is added (e.g., see [12] and references therein).  Furthermore, we
point out that the limiting of the energy spectrum at small scales due
to $\a$ does not, in itself, allow the flow to be resolved on a
coarser grid.  

The latter notwithstanding, we remark that the averaged Euler
equations are useful in better understanding the limit of inviscid
fluid flow, since the averaged Euler equations with viscosity, unlike
the Euler equations, converge regularly to the solutions of the
inviscid system\cite{Shk2}.  That is, for an arbitrary but fixed time
interval, we can choose $\a$ small enough so that the solution of the
averaged Euler equations are uniformly within any a priori chosen
error of the Euler equations\cite{OS} and then consider the zero
viscosity limit of the viscous, averaged Euler equations.

Finally, we note that as a theoretical model of fluid turbulence, the
averaged Euler equations in two dimensions possess other remarkable
features.  For example, for any initial condition and fixed time
interval, one may choose the number of modes $k_{max}$ large enough so
as to be arbitrarily close to the exact solution of the averaged Euler
equations without the addition of viscosity.  For such large
$k_{max}$, and in simulations of unforced decaying turbulence, the
averaged Euler equations exhibit a fundametal feature of
two-dimensional turbulence: a sharp decrease in enstrophy $Z$ during
the first few large eddy turnover times.  This is extremely
interesting, because, while it is necessary to add viscosity to the
Euler equations to obtain similar behavior, the averaged Euler
equations can reproduce this behavior while exactly conserving an
energy.  We shall report further on such inviscid simulations in
future publications.

The authors thank Jerry Marsden, and David Montgomery for 
extensive discussions on a variety of issues related to this article. 
BTN was supported by the Climate Change Prediction Program of the Department
of Energy, and SS was partially supported by NSF-KDI grant ATM-98-73133
and the Alfred P. Sloan Research Fellowship.

\newpage

\end{document}